\documentclass[12pt,a4paper]{amsart}
\usepackage[utf8]{inputenc}
\usepackage{a4wide,amssymb,xcolor,array,hyperref}
\usepackage{amsmath}
\usepackage{float}
\usepackage{geometry,graphicx,xcolor}
\usepackage{amssymb,color,esint,epic, graphicx,tikz}
\usepackage{parskip}
\usepackage{tcolorbox}
\usepackage{hyperref}
\usepackage{amsthm, amsfonts,enumerate}

\hypersetup{
 colorlinks,
 linkcolor={blue!90!black},
 citecolor={red!80!black},
 urlcolor={blue!50!black}
}

\newcommand{\vertiii}[1]{{\left\vert\kern-0.25ex\left\vert\kern-0.25ex\left\vert #1 
    \right\vert\kern-0.25ex\right\vert\kern-0.25ex\right\vert}}

\theoremstyle{plain}
\newtheorem*{thm}{Theorem}

\begin{document}

\title{An Improved Lower Bound for the Logarithmic Energy on $\mathbb S^2$}

\author{J. Marzo}
\address{Dept.\ Matem\`atica i Inform\`atica,
	Universitat  de Barcelona and BGSMath,
	Gran Via 585, 08007 Bar\-ce\-lo\-na, Spain}
\email{\href{mailto:jmarzo@ub.edu}{\texttt{jmarzo@ub.edu}}}

\thanks{The authors have
	been partially supported by grant MTM2017-83499-P by the 
	Ministerio de Econom\'{\i}a y Competitividad, Gobierno de Espa\~na and by the 
	Generalitat de Catalunya (project 2017 SGR 358).}

\begin{abstract}
In this short note, we employ well-known results to improve the lower bound for the constant associated with the linear term in the asymptotic expansion of the minimal logarithmic energy on the sphere.
\end{abstract}

\date{\today}

\maketitle

Given \( N \geq 2 \), the minimal logarithmic energy is defined as
\[
\mathcal{E}_N = \min_{x_1, \dots, x_N \in \mathbb{S}^2} E(x_1, \dots, x_N),
\]
where
\[
E(x_1, \dots, x_N) = \sum_{i \neq j} \log \frac{1}{|x_i - x_j|}
\]
is the logarithmic energy of the configuration \( x_1, \dots, x_N \in \mathbb{S}^2 \); see \cite[Chapter 6]{BHS19}.

It was shown by Tsuji, in his proof of the elliptic version of the fundamental potential theory \cite[p.93]{Tsu59}, that
\[
\left( \max_{x_1, \dots, x_N \in \mathbb{S}^2} \prod_{i \neq j} |x_i - x_j| \right)^{\frac{1}{N(N-1)}} \searrow \frac{2}{\sqrt{e}}, \quad N \to +\infty,
\]
and therefore, taking logarithms,
\[
\lim_{N \to +\infty} \frac{\mathcal{E}_N}{N(N-1)} = \frac{1}{2} - \log 2.
\]

This result implies that
\[
\mathcal{E}_N = I(\sigma) N^2 + o(N^2), \quad N \to +\infty,
\]
where
\[
I(\sigma) = \int_{\mathbb{S}^2} \int_{\mathbb{S}^2} \log \frac{1}{|x - y|} \, d\sigma(x) \, d\sigma(y),
\]
and \( \sigma \) denotes the normalized surface measure on the sphere.

Interest in the asymptotic expansion of the minimal energy \( \mathcal{E}_N \) grew significantly following the formulation of Smale’s 7th problem~\cite{Sma00}, which asks for the construction of configurations \( x_1, \dots, x_N \in \mathbb{S}^2 \) such that
\[
E(x_1, \dots, x_N) - \mathcal{E}_N \leq C \log N,
\]
for some universal constant \( C > 0 \). The logarithmic term in this bound arises from the upper bound established in~\cite{SS93}:
\[
\sqrt{N(N+1)} \, e^{\frac{1}{2}(E(x_1, \dots, x_N) - \mathcal{E}_N)},
\]
which estimates the condition number of a polynomial whose roots are the stereographic projections of the points \( x_1, \dots, x_N \in \mathbb{S}^2 \). This estimate was motivated by the goal of constructing polynomials with polynomially bounded condition numbers.

A family of polynomials achieving optimal condition numbers was eventually constructed in~\cite{BEMOC21, BL22}; however, Smale’s 7th problem remains far from being solved—see~\cite[6.7]{BHS19} for further discussion.

Thanks to the work of many mathematicians \cite{Wag89, KS98, RSZ94, BHS12, BS16, Ste22, Lau21, SM76, LN75, BL23} (see also the monograph \cite[Chapter 6]{BHS19} for further details), the best available asymptotic expansion is:
\[
\mathcal{E}_N = I(\sigma) N^2 - \frac{1}{2} N \log N + C_{\log} N + o(N), \quad N \to +\infty,
\]
with
\[
-0.0568528\ldots = \log 2 - \frac{3}{4} \leq C_{\log} \leq C_{BHS},
\]
where
\[
C_{BHS} = 2 \log 2 + \frac{1}{2} \log \frac{2}{3} + 3 \log \left( \frac{\sqrt{\pi}}{\Gamma(1/3)} \right) = -0.0556053\ldots
\]

It was conjectured in \cite{BHS12} that \( C_{\log} = C_{BHS} \). An equivalent conjecture was independently formulated in \cite{SS12}; see \cite{BS16} for the equivalence.

The existence of the constant \( C_{\log} \) and the upper bound were established in \cite{BS16}. The lower bound was first obtained in \cite{Lau21}, who adapted earlier results from \cite{SM76} and \cite{LN75} to the spherical setting. A direct proof of this lower bound, still based on \cite{LN75}, was later given in \cite{BL23}. Our result improves the bound from \cite{Lau21,BL23}.

\begin{thm}
There exists a constant \( C_{\log} > 0 \) such that
\[
\mathcal{E}_N = I(\sigma) N^2 - \frac{1}{2} N \log N + C_{\log} N + o(N), \quad N \to +\infty,
\]
and
\[
-0.0568456\ldots = \tilde{C} \leq C_{\log} \leq C_{BHS} = -0.0556053\ldots,
\]
where
\[
\tilde{C} = \log 2 - \frac{3}{4} + \frac{1}{162} \left( \sqrt[4]{3} \sqrt{2\pi} \left( 2 + 3 \tanh^{-1}\left( \frac{1}{2} \right) \right) - 12 \right)^2.
\]
\end{thm}

\section{Proof}

First we recall the  proof in \cite{BL23}. The first ingredient is an energy decomposition of the smeared out measure. 
Denote 
for $x_1,\dots , x_N\in \mathbb S^2$ and $\epsilon>0$ the signed measure
$$\mu=\mu_{\{x_i\},\epsilon}=\frac{1}{N}\sum_{i=1}^N \mu_i - d\sigma,\;\;\mbox{for}\;\; \mu_i=\frac{\chi_{B(x_i,\epsilon/\sqrt{N})}}{\sigma(B(x_i,\epsilon/\sqrt{N}))}d\sigma,$$
where $B_i=B(x_i,\epsilon/\sqrt{N})$ is the spherical cap with center $x_i$ and radius $\epsilon/\sqrt{N}.$
In what follows we denote $G(x,y)=\log \frac{1}{|x-y|}.$  
Then for the energy 
$$I(\mu)=\int_{\mathbb S^2} \int_{\mathbb S^2} G(x,y) d\mu(x) d\mu(y)$$ we have, as $\mu(\mathbb S^2)=0,$ that
$$N^2 I(\mu)=
\sum_{i\neq j}\int_{\mathbb S^2}\int_{\mathbb S^2}G d\mu_i d\mu_j+\sum_{i}\int_{\mathbb S^2}\int_{\mathbb S^2}G d\mu_i d\mu_i-N\sum_{i}\int_{\mathbb S^2}\int_{\mathbb S^2}G d\mu_i d\sigma.$$

From the computations in \cite{BL23} we know that for $x,y\in \mathbb S^2$ and $a>0$
$$\frac{1}{\sigma(B(x,a))^2}\int_{B(x,a)}\int_{B(x,a)}G d\sigma d\sigma= 
-\kappa+\ln \sin \frac{a}{2}  +\cot^2 \frac{a}{2}\left(\frac{1}{2}+\cot^2 \frac{a}{2}\ln \cos \frac{a}{2}\right),$$
and 
$$\frac{1}{\sigma(B(x,a))^2}\int_{B(x,a)}\int_{B(y,a)}G d\sigma d\sigma\le G(x,y)+ 
1+2\cot^2 \frac{a}{2}\ln \cos \frac{a}{2}$$
with equality if $B(x,a)$ and $B(y,a)$ are disjoint.

Therefore by taking $a=\epsilon/\sqrt{N}$ and considering the Taylor expansion of the functions we get, for all $\epsilon>0,$
\begin{equation}\label{decomp}
\sum_{i\neq j}\log \frac{1}{|x_i-x_j|}\ge \kappa N^2-\frac{N}{2}\log N+I(\mu) N^2
+\left(-\frac{1}{4}+\log \epsilon-\frac{\epsilon^2}{8}\right) N +o(N).
\end{equation}

The lower bound by Lauritsen, as it is proved in \cite{BL23}, follows then by using that $I(\mu)\ge 0,$ see \cite{Ski15,BCC19}, and using that the function 
$$u(\epsilon)=-\frac{1}{4}+\log \epsilon-\frac{\epsilon^2}{8}$$ 
has its maximum at $\epsilon=2.$

Our goal is now to get a better bound for the constant of the linear term by proving a non-trivial lower bound, of the right order, for the energy $I(\mu )$ in (\ref{decomp}).

Recall that for probability measures $\alpha,\beta$ on the sphere
the Wasserstein distance is defined by
$$W_1(\alpha,\beta)=
\inf_{\Pi \in \mbox{Coup}(\alpha,\beta)} \int_{\mathbb S^2\times \mathbb S^2} d(x,y) \, d\Pi (x,y)$$
where $\mbox{Coup}(\alpha,\beta)$ are the transport plans $\Pi$ from $\alpha$ to $\beta$ i.e. $\Pi(\mathbb S^2\times \mathbb S^2)=1$ and 
$\Pi(A\times \mathbb S^2)=\alpha(A),$ $\Pi(\mathbb S^2\times A)=\beta(A)$ for all borel subsets $A\subset \mathbb S^2,$ \cite{Vil03}. 

Now as a particular instance of \cite[Lemma 3.2]{GZ19} we get that 
$$W_1( \frac{1}{N} \sum_{j=1}^N \frac{\chi_{B_i}}{\sigma(B_i)} d\sigma,d\sigma)^2\le 2 \int_{\mathbb S^2}\int_{\mathbb S^2}
\log \frac{1}{|x-y|}d\mu (x)d\mu (y)=2I(\mu).$$
Observe that we use a different normalization of the Green function and the space than the ones in \cite{GZ19}.

By the Kantorovich-Rubinstein theorem \cite[Th. 1.14]{Vil03}
$$W_1( \alpha, \beta)=\sup \left\{ \left| \int_{\mathbb S^2}f d\alpha-  \int_{\mathbb S^2}f d\beta \right|  \;:\;  \|  f \|_{\mbox{Lip}}=\sup_{x\neq y}\frac{|f(x)-f(y)|}{d(x,y)}\le 1 \right\}$$
where $d(\cdot ,\cdot )$ is the geodesic distance.

Therefore for the Lipschitz function (with Lipschitz constant 1)
$$f(x)=d(x,\cup_{i=1}^N B_i)$$ we get that
$$W_1( \frac{1}{N} \sum_{j=1}^N \frac{\chi_{B_i}}{\sigma(B_i)} d\sigma,d\sigma)\ge \int_{\mathbb S^2}f(x)d\sigma(x).$$

To bound the integral above we use the following classical result by Fejes T\'oth (see \cite[5.8]{FTFTK23}).

\begin{thm}{(Fejes T\'oth \cite[Th. 2]{FT48})}
Let $\Phi(s)$ be an increasing function of $0\le s\le \pi/2.$ Let $x_1,\dots , x_N\in \mathbb S^2$ not all of them on an hemisphere. Then
$$\int_{\mathbb S^2} \Phi(d(x,\{x_i \}_{i=1}^N))dx \ge 
(2N-4)\int_{T} \Phi(d(x,\{a,b,c \}))dx $$
where $T$ is an equilateral spherical triangle with vertices $a,b,c\in \mathbb S^2$ and area $2\pi/(N-2).$ 
\end{thm}

Indeed, take
$$f(x)=\Phi(d(x,\{x_i \}_{i=1}^N))$$
for
$\Phi(s)=\max\{ s-\epsilon/\sqrt{N},0 \}$
and then
$$ \hskip -0.2cm N^2 G(\mu)\ge \frac{N^2}{2} W_1(\mu+d\sigma,d\sigma)^2
\ge  \frac{N^2}{2}\left( (2N-4) \int_{T} \Phi(d(x,\{a,b,c \}))d\sigma(x) \right)^2.
$$
where $T$ is an equilateral spherical triangle of area $2\pi/(N-2).$

What it is left now is to compute the integral in the above equilateral spherical triangle.

In terms of the interior angles $A$ of the triangle $T$ we have
$$\frac{2\pi}{N-2}=3A-\pi$$
and then
$$A=\frac{1}{3}\left( \frac{2\pi}{N-2}+\pi\right).$$
L'Huilier's formula relates the area of $T$ with the side lengths $\alpha$ (angles) giving
$$\tan \frac{\pi}{2(N-2)}=\left( \tan \frac{3\alpha}{4}\tan^3 \frac{\alpha}{4}\right)^{1/2},$$
and we get 
$$\alpha= \sqrt{\frac{8\pi}{\sqrt{3}}}\frac{1}{\sqrt{N}}+o(N^{-1/2}).$$

To compute the integral 
$$\int_{T} \Phi(d(x,\{a,b,c \}))d\sigma(x)$$
we split the triangle in Voronoi cells $T=\mathcal{T}_a\cup \mathcal{T}_b\cup \mathcal{T}_c$ where 
$$\mathcal{T}_a=\{ x\in T \;:\; d(x,a)<\min\{ d(x,b),d(x,c) \}  \}.$$

By using normal coordinates centered at $a$ we have
$$\int_{\mathcal T_a} \Phi(d(x,a)) d\sigma(x)=2
\int_{\mathcal T_a^+} \Phi(d(x,a)) d\sigma(x)=
\frac{1}{2\pi}\int_0^{A/2}\int_{0}^{h_\theta}r \Phi(r)dr d\theta.$$
By the spherical law of cosines we get $\tan h_\theta=\tan(\alpha/2)/\cos \theta$ and then $h_\theta\sim \alpha/2\cos \theta.$

Putting all together 
$$\int_{T} \Phi(d(x,\{a,b,c \}))d\sigma(x)=\frac{3}{2\pi}
\int_0^{\pi/6}\int_{\epsilon/\sqrt{N}}^{\alpha/2\cos \theta}r(r-\frac{\epsilon}{\sqrt{N}})drd\theta+o(N^{-3/2}).$$
For 
$C=\sqrt{\frac{2\pi}{\sqrt{3}}},$ then the integral above equals 
$$\frac{3}{2\pi  N^{3/2}}
\int_0^{\pi/6}\int_{\epsilon}^{C/\cos \theta}s(s-\epsilon)ds d\theta=
\frac{3}{2\pi  N^{3/2}}
\int_0^{\pi/6}\left( \frac{C^3}{3\cos^3 \theta}-\frac{C^2 \epsilon}{2\cos^2 \theta}+\frac{\epsilon^3}{6}  \right) d\theta$$
$$=\frac{3}{2\pi  N^{3/2}}  \left[
\frac{C^3}{3} \int_0^{\pi/6} \frac{1}{\cos^3 \theta} d\theta - \frac{C^2 \epsilon}{2} \int_0^{\pi/6} \frac{1}{\cos^2 \theta} d\theta+\frac{\pi \epsilon^3}{36}  \right].$$

For $\epsilon<\sqrt{\frac{8\pi}{3\sqrt{3}}}=C/\cos(\pi/6)=2.19927\dots $
(for bigger $\epsilon$ the function $ \Phi(d(x,\{a,b,c \}))$ is constant zero)
$$N^2 \int_{\mathbb S^2}\int_{\mathbb S^2}
\log \frac{1}{|x-y|}d\mu_{\epsilon,N}(x)d\mu_{\epsilon,N}(y)\ge
\frac{9 N}{2\pi^2}\left[
\frac{C^3}{3} C_1 - \frac{C^2 \epsilon}{2} C_2+\frac{\pi \epsilon^3}{36}  \right]^2=v(\epsilon),$$
with 
$$C_1=\int_0^{\pi/6} \frac{1}{\cos^3 \theta} d\theta=\frac{2+3\tanh^{-1}(1/2)}{6},\;\;C_2=\int_0^{\pi/6} \frac{1}{\cos^2 \theta} d\theta=\frac{1}{\sqrt{3}}.$$

\begin{figure}[htbp]
\centering
\begin{minipage}[b]{0.48\textwidth}
  \centering
  \includegraphics[width=\textwidth]{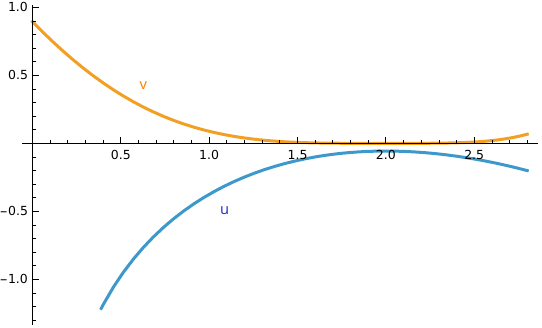}
  \caption{Functions $u$ and $v$}
  \label{fig1}
\end{minipage}
\hfill
\begin{minipage}[b]{0.48\textwidth}
  \centering
  \includegraphics[width=\textwidth]{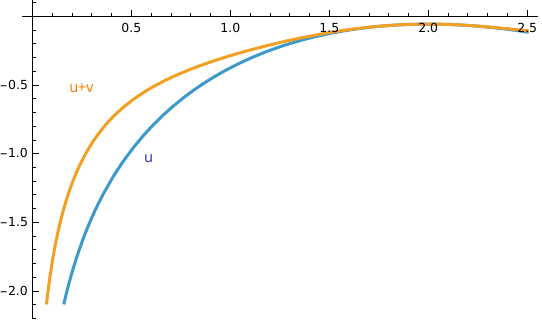}
  \caption{Functions $u$ and $u+v$}
  \label{fig2}
\end{minipage}
\end{figure}

Finally, the linear coefficient
$$\left(-\frac{1}{4}+\log \epsilon-\frac{\epsilon^2}{8}\right)+\frac{9 }{2\pi^2}\left[
\frac{C^3}{3} C_1 - \frac{C^2 \epsilon}{2} C_2+\frac{\pi \epsilon^3}{36}  \right]^2$$
has a maximum at $\epsilon=2$ (see Figures 1,2) and the value is
$$u(2)+v(2)=-0.0568456\dots>\log 2-\frac{3}{4}=u(2).$$\qed

{\bf Remarks} 1. It seems plausible that the lower bound could be further improved by considering the Wasserstein-2 distance. Indeed, one has the identity
$$\left\| \frac{1}{N} \sum_{j=1}^N \mu_i -d\sigma\right\|_{-1}^2= 2 I(\mu),$$
for 
$$\| \alpha \|_{-1}=\sup \left\{ \left| \int_{\mathbb S^2} f d\alpha \right| \;:\; f\in \mathcal{C}^\infty(\mathbb S^2) \;\mbox{and}\; \|\nabla f \|_{2}=\| f \|_{H^{1}}\le 1 \right\},$$
and by \cite{Pey18}
$$\inf_{\Pi \in \mbox{Coup}(\alpha,\beta)}\left(\int_{\mathbb S^2\times \mathbb S^2} d(x,y)^2 \, d\Pi (x,y) \right)^{1/2}=W_2( \alpha,\beta)\le 2 \| \alpha-\beta\|_{-1}.$$

The main challenge then becomes the analysis of the asymptotic behavior of a smeared quantization problem, namely
$$\lim_{N\to +\infty } \inf_{x_1,\dots, x_N\in \mathbb S^2}W_2(\frac{1}{N}\sum_{i}\frac{\chi_{B_i}}{\sigma(B_i)}d\sigma,d\sigma)N^{-1/2}.$$
See \cite{Klo12,GL00} for the standard quantization problem.

2. The fact that $I(\mu)$ is, for good configurations, of order $N^{-1}$ and it contributes to the constant of the linear term is already been observed in \cite{Wol}, see \cite{MM21}.

3. The process of smearing out the point masses can be done also in terms of other functions (not only characteristic functions of spherical caps). It was developed in \cite{Ste22} for the heat kernel.

\end{document}